\documentstyle{article}
\title{Intersection homology Betti numbers}
\author{Alan H. Durfee
    \thanks{Partially supported by NSF Grant DMS-8901903.
        Email: adurfee@mhc.mtholyoke.edu}
   \\ Department of Mathematics \\ Mount Holyoke College \\ South Hadley, MA
        01075
   }
\date{11 February 1993}

\def\endofproof{$\square$}
\def\euler{{\cal X}}

\input amssym.def
\input amssym.tex
\newtheorem{assertion}{Assertion}

\newtheorem{proposition}[assertion]{Proposition}

\newtheorem{lemma}[assertion]{Lemma}

\newcommand{\assertionnumber}{ \addtocounter{assertion}{1}
\thesection.\arabic{assertion}   }

\newenvironment{proof}{\bigskip \noindent {\bf Proof.}}{\endofproof \bigskip}
\begin{document}
\maketitle

%%%%%%%%%%%%%%%%%%%%%%%%%%%%%%%%%%%%%%%

\begin{abstract}
A generalization of the formula of Fine and Rao for the ranks of the
intersection homology groups of a complex algebraic variety is given.  The
proof uses geometric properties of intersection homology and mixed Hodge
theory.
\end{abstract}

\bigskip

%%%%%%%%%%%%%%%%%%%%%%%%%%%%%%%%%%%%%%%%%%

The middle-perversity intersection homology with integral coefficients of a
compact complex  n-dimensional algebraic variety $X$ with isolated
singularities is well known to be \cite{GM1}
\begin{equation}
\label{iso-ih-computation}
IH_k(X) \cong  \left\{ \begin{array}{ll}
               H_k(X)                     & \mbox{for $k > n$}  \\
               Im \,\{ H_n(U) \to H_n(X) \} & \mbox{for $k = n$} \\
                H_k(U)                     & \mbox{for $k < n$}
             \end{array}
           \right.
\end{equation}
where $U$ is the complement of the singular set in $X$.
Fine and Rao \cite{FR} find a formula for the ranks of the intersection
homology groups of $X$ in terms of the ordinary Betti numbers of a resolution
of singularities of $X$.
The formula is an application of (\ref{iso-ih-computation}), which is a purely
geometric fact, and mixed Hodge theory.

This paper contains a mild generalization of their formula to the case when the
variety has non-isolated singularities.
First, (\ref{iso-ih-computation}) is generalized to nonisolated singularities
(Proposition \ref{ih-computation}).
The proof uses a basic geometric property of intersection homology (Lemma
\ref{cohomology-vanishing}) and a standard exact sequence of mixed Hodge theory
(Lemma \ref{mhs-computation}).  A generalization of the Fine-Rao formula for
the Betti numbers (Proposition \ref{betti-formula}) follows easily. The proof
of the formula given in this paper makes clear the separate roles played by
geometry and mixed Hodge theory.

%%%%%%%%%%%%%%%%%%%%%%%%%%%%%%%%%%%%

Let $X$ be a complex algebraic variety which is not necessarily compact (i.e.,
a reduced separable scheme of finite type over the complex numbers), and
suppose that the complex dimension of $X$ is $n$.  Let $Z$ be a closed
subvariety of dimension $\leq d$ and let $U = X \backslash Z$.
We let $ IH^k_c(X) = IH_{2n-k}(X)$ (respectively
$IH^k(X) =  IH^{BM}_{2n-k}(X)$), where
$IH_k(X) $ (respectively $IH_k^{BM}(X)$) denote the  middle-perversity
intersection homology group of $X$  (respectively middle-perversity
Borel-Moore intersection homology group),
and similarily for open subsets of $X$.
For compact $X$ with isolated singularitites, formula
(\ref{iso-ih-computation}) thus becomes
\[
IH^k(X) \cong  \left\{ \begin{array}{ll}
               H^k(X)                     & \mbox{for $k > n$}  \\
               Im \,\{ H^n(X) \to H^n(U) \} & \mbox{for $k = n$} \\
                H^k(U)                     & \mbox{for $k < n$}
             \end{array}
           \right.
\]

Let $ IC^{top}_X $ be the sheaf complex defined by
$(IC^{top}_X)^k (V) = IC^{BM}_{-k} (V) $
for $V \subset X$ open, so that $IH^k(X)$ is the hypercohomology group
$H^k(X, IC^{top}_X)$, etc.
Let
$i: Z \hookrightarrow X \ \ \mbox{and} \ \ j: U \hookrightarrow X$.
We abbreviate $Ri_*$ by writing $i_*$ and so forth, and  we let
$
IC_X = IC^{top}_X [-n]
$as in \cite{FP}.

The following lemma expresses a basic geometric property of intersection
homology.
If $Z$ has a mapping cylinder neighborhood $N$ in $X$ (see \cite{DS}), then a
geometric argument similar to \cite[1.7]{GM.Morsetheory} shows the equivalent
result
$IH_k^{BM}(N) = 0$ for $k \leq n-d$ and
$IH_k(N) = 0$ for $k \geq n+d$.

\begin{lemma}
\label{cohomology-vanishing}
If $X$ is an algebraic variety of dimension $n$ and $U$ an open subvariety with
$Z = X \backslash U$ a closed subvariety of dimension $\leq d$,
then (with integer coefficients)
\begin{eqnarray*}
IH_c^k(X)  &  \cong  & IH_c^k(U) \ \ \mbox{for $k > n+d$}  \\
IH_c^{n+d}(X)  &  \twoheadleftarrow  &  IH_c^{n+d}(U)   \\
IH^{n-d}(X)   &  \hookrightarrow  & IH^{n-d}(U)   \\
IH^k(X)   &  \cong  &  IH^k(U) \ \ \mbox{for $k < n-d$}
\end{eqnarray*}
\end{lemma}

\begin{proof}
There is an exact triangle \cite[1.4.3.4]{FP}
\begin{equation}
\label{eqn1}
j_!j^*IC_X \to IC_X \to i_* i^* IC_X
\end{equation}
Taking $H^{k-n}_c (X, -)$ gives
$$
\to H^{k-n}_c(X, j_!j^* IC_X) \to H^{k-n}_c (X, IC_X) \to H^{k-n}_c (X, i_*i^*
IC_X) \to
$$
which gives
\begin{equation}
\label{eqn2}
\to IH^k_c(U) \to IH^k_c(X) \to H^{k-n}_c(Z, i^* IC_X) \to
\end{equation}
Now $H^{k-n}(i^*IC_X)| Z = 0$ as sheaf for $k-n \geq -d$ by \cite[p. 9]{FP}.
Since $Z$ has complex dimension $\leq d$, we have $H^{k-n}(Z, i^*IC_X) = 0$ for
$k-n \geq d$.
This proves the first part of the lemma.

Similarly there is an exact triangle
\begin{equation}
\label{eqn3}
i_*i^!IC_X \to IC_X \to j_* j^! IC_X
\end{equation}
Taking $H^{k-n} (X, -)$ gives
$$
\to H^{k-n}(X, i_*i^! IC_X) \to H^{k-n}(X, IC_X) \to H^{k-n} (X, j_*j^! IC_X)
\to
$$
which gives
\begin{equation}
\label{eqn4}
\to H^{k-n}(Z, i^! IC_X) \to IH^k(X) \to IH^k(U) \to
\end{equation}
Now $H^{k-n}(i^!IC_X)| Z = 0$ as sheaf for $k-n \leq -d$ by \cite[p. 9]{FP}, so
$H^{k-n}(Z, i^!IC_X) = 0$ for $k-n \leq -d$.
This proves the last part of the lemma.
\end{proof}

%%%%%%%%%%%%%%%%%%%%%%%%%%%%%%%%%%%%%%%%%

Next we need a basic property of the mixed Hodge structure on intersection
cohomology.
By work of Saito \cite{Saito.ICM}, the intersection cohomology groups of an
algebraic variety have a mixed Hodge structure.
In particular they have a weight filtration
$
\dots \subset W_m \subset W_{m+1} \subset \dots
$
This filtration is functorial for algebraic maps, and
$Gr^W_m = W_m / W_{m-1}$ is an exact functor.

\begin{lemma}
\label{mhs-computation}
If $X$ is a complex algebraic variety, and $U$ an open subvariety, then for all
$k$ (with rational coefficients):
\begin{enumerate}
\item $Gr^W_k IH^k_c (U) \hookrightarrow Gr^W_k IH^k_c (X)$
\item $Gr^W_k IH^k (X) \twoheadrightarrow Gr^W_k IH^k(U)$
\end{enumerate}
\end{lemma}

The first part of this lemma is a generalization of the fact that for a smooth
compactification $X$ of a smooth $U$, the classes in $H^k(U)$ of weight $k$ are
restrictions of classes in $H^k(X)$ of weight $k$
\cite[3.2.17]{Deligne.Hodge-II}.
The proof of the lemma uses Saito's theory of mixed Hodge modules
\cite{Saito.ICM}, which is parallel to the theory of mixed perverse sheaves
\cite{FP}.  We recall that $IC_X$ as mixed Hodge module has pure weight $n$,
where $n$ is the dimension of $X$.  Also if $f$ is a map of algebraic
varieties, then $f_!$ and $f^*$ preserve weight $\leq m$, and $f_*$ and $f^!$
preserve weight $\geq m$.
Finally, if $K$ is a complex of weight $\leq m$ (respectively $\geq m$), then
$H^k(K)$ has weight $\leq m+k$ (respectively $\geq m+k$).

\begin{proof}
Let
$
a: X \to pt
$.
The exact triangle (\ref{eqn1}) is an exact triangle of mixed Hodge modules
\cite[4.4.1]{Saito.MHM}.
Taking $a_!$ and $H^{k-n}( - )$ is the same as taking $H^{k-n}_c (X, -)$,
so all terms of exact sequence (\ref{eqn2})
have weight $\leq k$, and taking $Gr^W_k$ proves the first part of the lemma.

Similarly taking $a_*$ and $H^{k-n}(-)$ of exact triangle (\ref{eqn3}) is the
same as taking $H^k (X, -)$,
so all terms of exact sequence (\ref{eqn4})
have weight $\geq k$, and taking $Gr^W_k$ proves the second part of the lemma.
\end{proof}

Combining Lemmas \ref{cohomology-vanishing} and \ref{mhs-computation} yields
the following proposition.

\begin{proposition}
\label{ih-computation}
If $X$ is a compact algebraic variety of dimension $n$ and $U$ an open smooth
subvariety with $Z = X \backslash U$ of dimension $\leq d$, then (with rational
coefficients)
$$
IH^k(X) \cong  \left\{ \begin{array}{ll}
               H^k_c(U)                     & \mbox{for $k > n+d$}  \\
                                           &                       \\
               Gr^W_{n+d} H^{n+d}_c(U)      & \mbox{for $k = n+d$} \\
					  &                      \\
               Gr^W_{n-d} H^{n-d}(U)      & \mbox{for $k = n-d$} \\
                                           &                      \\
               H^k(U)                     & \mbox{for $k < n-d$}
             \end{array}
           \right.
$$
\end{proposition}

The Proposition follows since
$IH^k(X)$ has pure weight $k$.  Note that $H^k_c(U)$ can be replaced by
$H^k(X,Z) \cong H^k(X)$ since $k \geq n+d$.
The  proposition depends on mixed Hodge theory and hence is valid only over the
rational numbers.
Is there an integral version of this proposition, with a geometric proof
similar to Formula \ref{iso-ih-computation}?

%%%%%%%%%%%%%%%%%%%%%%%%%%%%%%%%%%%%%%%

The weighted Euler characteristic of a space $Y$ with mixed Hodge structure is
defined by
$$
\euler_m(Y) = \sum _{i} (-1)^i dim \, Gr_m^W H^i(Y)
$$
where $W$ is the weight filtration.
Similarly we define
$\euler_m^c(Y)$ using $H^i_c(Y)$ in place of $H^i(Y)$.
The following proposition gives a formula for the ranks of the outer
intersection homology groups of an algebraic variety in terms of the ordinary
Betti numbers of a resolution of its singularities.

\begin{proposition}
\label{betti-formula}
Let $X$ be a compact variety of dimension $n$ and let $Z$ be a closed
subvariety of dimension $\leq d$ containing the singular set of $X$.
Let $\pi : \tilde X \to X$ be an algebraic map with $\tilde X$ smooth and
$\pi^{-1} (Z) = E = E_1 \cup \dots  \cup E_r $
a divisor with normal crossings.
Suppose that $\pi$ is an analytic isomorphism of $\tilde U = \tilde X
\backslash E$ to $U = X \backslash Z$.
Then for $m \geq n+d$
$$
dim \, IH^m(X) = dim \, H^m(\tilde X) - \sum_i (-1)^i dim \, H^m(E_i)
$$
\end{proposition}

\begin{proof}
By Proposition \ref{ih-computation} and the fact that $IH^k(X)$ is pure and the
weights on $H^k_c(U)$ are $\leq k$ for all $k$,
$$
(-1)^m dim \, IH^m(X) = (-1)^m dim \, Gr^W_mIH^m_c(U) = \euler^c_m(U)
$$
for $m \geq n+d$.  The exact sequence of mixed Hodge structures
$$
\to H^k_c(\tilde U) \to H^k(\tilde X) \to H^k(E) \to
$$
gives
$$
\euler^c_m(\tilde U) = \euler_m(\tilde X) - \euler_m(E)
$$
Since  $\tilde X$ is compact and smooth, its cohomology is pure and
$$
\euler_m(\tilde X) = (-1)^m \, dim \, H^m(\tilde X)
$$
The exceptional set $\tilde Z$ has a simplicial resolution
$$
E_\bullet =  \dots  \to E_{ [l] } \to E_{ [l-1] } \to \dots
$$
where
$$
E_{ [l] } = \coprod_{|I| = l+1} ( \bigcap_{i \in I} E_i )
$$
Thus
$$
\euler_m(\tilde Z) = \sum (-1)^i \euler_m(E_i) = \sum (-1)^i (-1)^m dim \, H^m
(E_i)
$$
which finishes the proof.
\end{proof}

By duality this gives the ranks of the groups in dimensions  $ \leq n-d$.
The same argument also gives a formula for the Hodge numbers.

What about the missing middle-dimensional Betti numbers?
The following is a rather naive estimate of the rank of intersection homology
groups.

\begin{proposition}
Let $X$ be a compact algebraic variety, let $U \subset X$ be a smooth open
subset, and let $\tilde X \to X$ be a resolution of singularities.  Then for
all $k$,
$$
max \{ dim \, Gr_k^W H^k(U) , \ dim \, Gr_k^W H^k_c(U) \}
\leq dim \, IH^k(X)
\leq dim \, H^k( \tilde X)
$$
\end{proposition}

This follows from Lemma \ref{mhs-computation} and the fact that the
decomposition theorem implies that the intersection homology of $X$ is a direct
summand of the ordinary homology of $ \tilde X$ \cite{FP}.
The lower bound is realized when the computation of intersection homology is
essentially reduced to (\ref{iso-ih-computation}) or Proposition
\ref{betti-formula}.
The upper bound is realized if the resolution $\tilde X \to X$ is small.  For
example, Schubert varieties have small resolutions \cite{BFL}.  If the
normalization of a variety is smooth then it is a small resolution as well.

%%%%%%%%%%%%%%%%%%%%%%%%%%%%%%%%%%%%%%%%%%%%%%

\bibliographystyle{alpha}

\end{document}